\documentclass[a4paper,oneside,11pt]{article}
\usepackage[utf8]{inputenc}
\usepackage[textwidth=14cm]{remipeyre}

\DeclareMathOperator\dist{\mathit{dist}}

\DeclareMathOperator\closure{cl}
\newcommand\id{\mathord{\mathit{Id}}}
\newcommand\pushfwd{\mathbin{{}_*}}
\newcommand\sphere[1]{\mathbb{S}^{#1}}
\newcommand\dHmo{\dot{H}^{-1}}
\newcommand\csf{\mathsf{c}}
\newcommand\Lipschitz{\mathit{Lip}}

\newcommand\stepone{\text{\ding{172}}}
\newcommand\steptwo{\text{\ding{173}}}
\newcommand\stepthree{\text{\ding{174}}}

\theoremstyle{remark}
\newtheorem*{acknowledgement}{Acknowledgement}

\title[$\smash{W_2}$ distance and $\smash{\dHmo}$ norm]{Comparison between $W_2$ distance and $\dHmo$ norm,\penalty-1
\ and localisation of Wasserstein distance}
\author{Rémi Peyre}
\date{\today}

\begin{document}

\maketitle

\begin{abstract}
It is well known that the quadratic Wasserstein distance $W_2 (\placeholder, \placeholder)$ is formally equivalent, for infinitesimally small perturbations, to some weighted $H^{-1}$ homogeneous Sobolev norm. In this article I show that this equivalence can be integrated to get non-asymptotic comparison results between these distances. Then I give an application of these results to prove that the $W_2$ distance exhibits some localisation phenomenon: if $\mu$ and $\nu$ are measures on $\RR^n$ and $\phi \colon \RR^n \to \RR_+$ is some bump function with compact support, then under mild hypotheses, you can bound above the Wasserstein distance between $\phi \cdot \mu$ and $\phi \cdot \nu$ by an explicit multiple of $W_2 (\mu, \nu)$.
\end{abstract}

\section*{Foreword}

This article is divided into two sections, each of which having its own introduction. \S~\ref{sec:equivalence} deals with general results of comparison between Wasserstein distance and homogeneous Sobolev norm, while \S~\ref{sec:localisation} handles an application to localisation of $W_2$ distance.

\section{Non-asymptotic equivalence between $W_2$ distance and $\dHmo$ norm}\label{sec:equivalence}

\subsection{Introduction}

In all this section, $M$ denotes a connected Riemannian manifold endowed with its distance $\dist (\placeholder, \placeholder)$ and its Lebesgue measure $\lambda$. Let us give a few standard definitions which will be at the core of our work:
\begin{itemize}
\item For $\mu, \nu$ two positive measures on $M$, denoting by $\Pi (\mu, \nu)$ the set of (positive) measures on $M \times M$ whose respective marginals are $\mu$ and $\nu$, for $\pi \in \Pi (\mu, \nu)$ one defines
\[
I (\pi) \coloneqq \int_{M \times M} \dist (x, y)^2 \pi (\dx{x}, \dx{y})
\]
and then
\[\label{eqn:W2}
W_2 (\mu, \nu) \coloneqq \inf \setst{I (\pi)}{\pi \in \Pi (\mu, \nu)}^{1/2}
.\]
$W_2$ is a (possibly infinite) distance, called the \emph{quadratic Wasserstein distance} \citep[\S~7.1]{topics}. Note that this distance is finite only between measures having the same total mass.
\item On the other hand, for $\mu$ a (positive) measure on $M$, if $f$ is a $\mathcal{C}^1$ real function on $M$, one denotes
\[\label{eqn:dHo}
\norm{f}_{\dot{H}^1 (\mu)} \coloneqq \paren[\Big]{\int_M \abs{\nabla{f}(x)}^2 \dx{\mu} (x)}^{1/2}
,\]
which defines a semi-norm; for $\nu$ a signed measure on $M$, one then denotes
\[\label{eqn:dHmo}
\|\nu\|_{\dHmo (\mu)} \coloneqq \sup \setst{\abs{\inner{f}{\nu}}}{\norm{f}_{\dot{H}^1(\mu)} \leq 1}
,\]
which defines a (possibly infinite) norm, which we will call the \emph{$\dHmo (\mu)$ weighted homogeneous Sobolev norm}. Note that this norm is finite only for measures having zero total mass. In the case $\mu$ is the Lebesgue measure, we will merely write \quotes{$\dHmo$} for \quotes{$\dHmo (\lambda)$}.
\end{itemize}

The $W_2$ Wasserstein distance is an important object in analysis; but it is non-linear, which makes it harder to study. For infinitesimal perturbations however, the linearised behaviour of $W_2$ is well known: if $\mu$ is a positive measure on $M$ and $\dx{\mu}$ is an infinitesimally small perturbation of this measure,%
\footnote{Beware that here $\dx{\mu}$ denotes a small measure on $M$, not the value of $\mu$ on a small area.}
one has formally (see \citep[\S~7.6]{topics} or \cite[\S~7]{OttoVillani})
\[\label{eqn:BenamouBrenier.infinitesimal}
W_2 (\mu, \mu + \dx{\mu}) = \norm{\dx{\mu}}_{\dHmo (\mu)} + o (\dx{\mu})
.\]
More precisely, one has the following equality, known as the \emph{Benamou–Brenier formula} \citep[Prop.~1.1]{BenamouBrenier}: for two positive measures $\mu, \nu$ on $M$,
\[\label{eqn:BenamouBrenier}
W_2 (\mu, \nu) = \inf \setst[\Big]{\int_0^1 \norm{\dx{\mu_t}}_{\dHmo (\mu_t)}}{\mu_0 = \mu,\ \mu_1=\nu}
.\]
Then, a natural question is the following: are there \emph{non-asymptotic} comparisons between the $W_2$ distance and the $\dHmo$ norm? Concretely, we are looking for inequalities like
\[
C_{\mrm{a}} \norm{\mu - \nu}_{\dHmo (\mu)} \leq W_2 (\mu, \nu) \leq C_{\mrm{b}} \norm{\mu - \nu}_{\dHmo (\mu)}
\]
for constants $0 < C_{\mrm{a}} \leq C_{\mrm{b}} < \infty$, under mild assumptions on $\mu$ and $\nu$.

\subsection{Controlling $W_2$ by $\dHmo$}

\begin{theorem}\label{thm:W2<dHmo}
For any positive measures $\mu, \nu$ on $M$,
\[
W_2 (\mu, \nu) \leq 2 \norm{\mu - \nu}_{\dHmo (\mu)}
.\]
\end{theorem}

\begin{proof}
We suppose that $\norm{\mu - \nu}_{\dHmo (\mu)} < \infty$, otherwise there is nothing to prove. For $t \in [0, 1]$, let
\[
\mu_t \coloneqq (1 - t) \mu + t \nu
,\]
so that $\mu_0 = \mu$, $\mu_1 = \nu$ and $\dx{\mu_t} = (\mu - \nu) \dx{t}$. Then, by the Benamou–Brenier formula \eqref{eqn:BenamouBrenier}:
\[
W_2 (\mu, \nu) \leq \int_0^1 \norm{\mu - \nu}_{\dHmo (\mu_t)} \dx{t}
.\]

Now, we use the following key lemma, whose proof is postponed:
\begin{lemma}\label{lem:dHmo(mu')vsdHmo(mu)}
If $\mu, \mu'$ are two measures such that $\mu' \geq \rho \mu$ for some $\rho > 0$, then $\norm{\placeholder}_{\dHmo (\mu')} \leq \rho^{-1/2} \norm{\placeholder}_{\dHmo (\mu)}$.%
\footnote{Beware that here \squo{$\placeholder$} stands for a \emph{measure}, not for a function: otherwise the formula would be false.—When $f$ is a function, $\norm{f}_{\dHmo (\mu)}$ stands for the $\dHmo (\mu)$ norm of the measure having density $f$ w.r.t.\ $\mu$.}
\end{lemma}
\noindent Here obviously $\mu_t \geq (1 - t) \mu$, so
\[
W_2 (\mu, \nu) \leq \int_0^1 (1 - t)^{-1/2} \norm{\mu - \nu}_{\dHmo (\mu)} \dx{t} = 2 \norm{\mu - \nu}_{\dHmo (\mu)}
,\]
\textsc{qed}.
\end{proof}

\begin{corollary}
If $\mu \geq \rho \lambda$ for some $\rho > 0$, then
\[
W_2 (\mu, \nu) \leq 2 \rho^{-1/2} \norm{\mu - \nu}_{\dHmo}
.\]
\end{corollary}

\begin{proof}
Just use that $\norm{\placeholder}_{\dHmo (\mu)} \leq \rho^{-1/2} \norm{\placeholder}_{\dHmo}$ by Lemma~\ref{lem:dHmo(mu')vsdHmo(mu)}.
\end{proof}

\begin{proof}[\proofname\ of Lemma~\ref{lem:dHmo(mu')vsdHmo(mu)}]
Take $\mu' \geq \rho \mu$ and let $\nu$ be a signed measure on $M$ such that $\mu + \nu$ is positive; then $\mu' + \rho \nu$ is also positive. For $m$ a measure on $M$, we denote by $\mathit{diag} (m)$ the measure on $M \times M$ supported by the diagonal whose marginals (which are equal) are $m$, i.e.:
\[
\paren[\big]{\mathit{diag} (m)} (A \times B) \coloneqq m (A \cap B)
;\]
with that notation,
\[
\pi \in \Pi (\mu, \mu + \nu) \Rightarrow \rho \pi + \mathit{diag} (\mu' - \rho \mu) \in \Pi (\mu', \mu' + \rho \nu)
,\]
and
\[
I \paren[\big]{\rho \pi + \mathit{diag} (\mu' - \rho \mu)} = \rho I (\pi)
.\]
Therefore, taking infima,
\begin{multline}
W_2 (\mu', \mu' + \rho \nu)^2 = \inf \setst{I (\pi')}{\pi' \in \Gamma (\mu', \mu' + \rho \nu)} \\
\leq \inf \setst[\big]{I \paren[\big]{\rho \pi + \mathit{diag} (\mu' - \rho \mu)}}{\pi \in \Gamma (\mu, \mu + \nu)} \\
= \rho \inf \setst{I (\pi)}{\pi \in \Gamma (\mu, \mu + \nu)}
= \rho W_2 (\mu, \nu)^2
.\end{multline}
For infinitesimally small $\nu$, it follows by Equation~\eqref{eqn:BenamouBrenier.infinitesimal} that $\norm{\rho \nu}_{\dHmo (\mu')}^2 \leq \rho \norm{\nu}_{\dHmo (\mu)}^2$, hence $\norm{\nu}_{\dHmo (\mu')} \leq \rho^{-1/2} \norm{\nu}_{\dHmo (\mu)}$. This relation remains true even for non-infinitesimal $\nu$ by linearity, which ends the proof.
\end{proof}

\begin{remark}
Lemma \ref{lem:dHmo(mu')vsdHmo(mu)} could also be proved very quickly by using the definition \eqref{eqn:dHo}-\eqref{eqn:dHmo} of the $\dHmo (\mu)$ norm. The proof above, however, has the advantage that it does not need the precise expression of $\norm{\placeholder}_{\dHmo (\mu)}$, but only the fact that it is the linearised $W_2$ distance.
\end{remark}

\subsection{Controlling $\dHmo$ by $W_2$}

\begin{theorem}
Assume $M$ has nonnegative Ricci curvature. Then for any positive measures $\mu, \nu$ on $M$ such that $\mu \leq \rho_0 \lambda$ and $\nu \leq \rho_1 \lambda$,
\[\label{eqn:Hmo<W2}
\norm{\mu - \nu}_{\dHmo} \leq \frac{2 (\rho_0^{1/2} - \rho_1^{1/2})}{\ln (\rho_0 \div \rho_1)} W_2(\mu,\nu)
.\]
(For $\rho_1 = \rho_0$, the right-hand side of \eqref{eqn:Hmo<W2} is to be taken as $\rho_0^{1/2} W_2(\mu, \nu)$ by continuity).
\end{theorem}

\begin{remark}
For $M = \RR^n$ a similar result was already stated in \citep[Proposition~2.8]{Loeper:06}, with a different proof.
\end{remark}

\begin{proof}
Let $(\mu_t)_{0 \leq t \leq 1}$ be the displacement interpolation between $\mu$ and $\nu$ (cf.\ \citep[chap.~7]{OldAndNew}), which is such that $\mu_0 = \mu,\ \mu_1 = \nu$ and the infimum in \eqref{eqn:BenamouBrenier} is attained with $\norm{\dx{\mu_t}}_{\dHmo (\mu_t)} = W_2 (\mu, \nu) \dx{t}\ \forall t$. Since Ricci curvature is nonnegative, the Lott–Sturm–Villani theory tells us that, denoting by $\norm{\mu}_\infty$ the essential supremum of the density of $\mu$ w.r.t.\ $\lambda$, one has $\norm{\mu_t}_\infty \leq \norm{\mu_0}_\infty^{1 - t} \norm{\mu_1}_\infty^t = \rho_0^{1 - t} \rho_1^t$ (see \citep[Corollary~17.19]{OldAndNew} or \citep[Lemma~6.1]{CorderoerausquinAlii}); so that $\norm{\placeholder}_{\dHmo} \leq \rho_0^{(1 - t) / 2} \rho_1^{t / 2} \norm{\placeholder}_{\dHmo (\mu_t)}$ by Lemma~\ref{lem:dHmo(mu')vsdHmo(mu)}.

Then, by the integral triangle inequality for normed vector spaces,
\begin{multline}
\norm{\mu - \nu}_{\dHmo} = \norm[\Big]{\int_0^1 \dx{\mu_t}}_{\dHmo}
\leq \int_0^1 \norm{\dx{\mu_t}}_{\dHmo} \\
\leq \int_0^1 \rho_0^{(1 - t) /2} \rho_1^{t /2} \norm{\dx{\mu_t}}_{\dHmo (\mu_t)}
= \paren[\Big]{\int_0^1 \rho_0^{(1 - t) /2} \rho_1^{t /2} \dx{t}} W_2 (\mu, \nu) \\
= \frac{2 (\rho_0^{1/2} - \rho_1^{1/2})}{\ln (\rho_0 \div \rho_1)} W_2 (\mu, \nu)
,\end{multline}
\textsc{qed}.
\end{proof}

\begin{remark}
Taking into account the dimension $n$ of the manifold $M$, the bound on $\norm{\mu_t}_\infty$ could be refined into
\[
\norm{\mu_t}_\infty \leq \paren[\big]{(1 - t) \norm{\mu_0}_\infty^{-1 \div n} + t \norm{\mu_1}_\infty^{-1 \div n}}^{-n}
,\]
which would yield a slightly sharper bound in Equation~\eqref{eqn:Hmo<W2}, namely:
\begin{multline}\label{eqn:Hmo<W2.refined}
\norm{\mu - \nu}_{\dHmo} \leq \paren[\Big]{\int_0^1 \paren[\big]{(1 - t) \rho_0^{-1 \div n} + t \rho_1^{-1 \div n}}^{-n /2} \dx{t}} W_2 (\mu, \nu) \\
= \begin{cases}
\frac{\rho_0^{1/2 - 1 \div n} - \rho_1^{1/2 - 1 \div n}}{(n /2 - 1) (\rho_1^{-1 \div n} - \rho_0^{-1 \div n})} W_2 (\mu, \nu) & n \geq 2; \\
\frac{\log (\rho_1 \div \rho_0)}{2 (\rho_0^{-1/2} - \rho_1^{-1/2})} & n = 2. \end{cases}
\end{multline}
For $n = 1$ it turns out that one can let tend $\rho_1 \to \infty$ in \eqref{eqn:Hmo<W2.refined} without making the integral diverge; which leads to a much more powerful result:
\begin{theorem}
When $M$ is an interval of $\RR$, then under the sole assumption that $\mu \leq \rho_0 \lambda$, one has for all positive measures $\nu$ on $M$:
\[
\norm{\mu - \nu}_{\dHmo} \leq 2 \rho_0^{1/2} W_2 (\mu, \nu)
.\]
\end{theorem}
\end{remark}

\begin{remark}
For $n \geq 2$ there is no hope to get a bound valid for all $\nu$, because then it can occur that $W_2 (\mu, \nu) < \infty$ but $\norm{\mu - \nu}_{\dHmo} = \infty$: for instance, take $\mu$ to be the uniform measure on the $2$-dimensional sphere and $\nu$ a Dirac mass.
\end{remark}

\section{Application to localisation of Wasserstein distance}\label{sec:localisation}

\subsection{Introduction}

In all this section, we work in the Euclidian space $\RR^n$, whose norm is denoted by $\abs{\placeholder}$. $\dist (x, A) \coloneqq \inf \setst{\abs{x - y}}{y \in A}$ denotes the distance between a point $x$ and a set $A$; $\C{A}$ denotes the complement of $A$; $\lambda$ denotes the Lebesgue measure. We will use the following notation to handle measures:
\begin{itemize}
\item For $\mu$ a measure on $\RR^n$ and $f \colon \RR^n \to \RR^n$ a measurable map, $f \pushfwd \mu$ denotes the pushforward of $\mu$ by $f$, that is, $\paren[\big]{f \pushfwd \mu} (A) \coloneqq \mu (f^{-1} (A))$.
\item For $\mu$ a measure on $\RR^n$ and $\phi \colon \RR^n \to \RR_+$ a nonnegative measurable function, $\phi \cdot \mu$ denotes the measure such that $\dx{\paren[\big]{\phi \cdot \mu}} (x) \coloneqq \phi(x) \dx{\mu} (x)$.
\end{itemize}
We will also use the following norms on measures:
\begin{itemize}
\item $\norm{\mu}_{\dHmo (\nu)}$ has the same definition as in \S~\ref{sec:equivalence};
\item $\norm{\mu}_1 \coloneqq \int_{\RR^n} \abs{\dx{\mu} (x)}$ is the total variation norm of $\mu$;%
\footnote{Note that in the case $\mu$ is a positive measure on $\RR^n$, then $\norm{\mu}_1$ is noting but $\mu (\RR^n)$.}
\item For $\mu$ a measure supported by $A$, we define
\[
\norm{\mu}_{L^2 (A)} \coloneqq \paren[\bigg]{\int_A \paren[\bigg]{\frac{\dx{\mu}}{\dx{\lambda}} (x)}^2 \dx{\lambda} (x)}^{1/2}
.\]
\end{itemize}

The goal of this section is to give an application of Theorem~\ref{thm:W2<dHmo} to the problem of \emph{localisation} of the quadratic Wasserstein distance. Morally, the question is the following: take two measures $\mu, \nu$ on $\RR^n$ being close to each other in the sense of $W_2$ distance; is it true that $\mu$ and $\nu$ remain close when you consider their restrictions to a subset of $\RR^n$? Concretely, if $\phi$ is a non-negative real function on $\RR^n$ with compact support (plus some technical assumptions to be specified later), we want to bound above $W_2 (a \phi \cdot \mu, \phi \cdot \nu)$ by some multiple of $W_2 (\mu, \nu)$—where, in the former expression, $a$ is a constant factor ensuring that $a \phi \cdot \mu$ and $\phi \cdot \nu$ have the same mass (for otherwise the distance between $\phi \cdot \mu$ and $\phi \cdot \nu$ is generically infinite).

This question, which was my initial motivation for the results of \S~\ref{sec:equivalence}, was asked to me by Xavier \textsc{Tolsa}, who needed such a result for his paper \citep{Tolsa.published} on characterizing uniform rectifiability in terms of mass transport. Actually Xavier managed to devise a proof of his own \citep[Theorem~1.1]{Tolsa.published}, but it was quite long (about thirty pages) and involved arguments of multi-scale analysis. With Theorem~\ref{thm:W2<dHmo} at hand, however, the reasoning becomes far more direct; moreover we will be able to relax some of the assumptions of Xavier's theorem.

\subsection{Statement of the theorem}

\begin{theorem}\label{thm:localisation.sphere}
Let $\mu, \nu$ be (positive) measures on $\RR^n$ having the same total mass; let $B$ be a ball of $\RR^n$ (whose radius will be denoted by $R$ when needed). Assume that on $B$, the density of $\mu$ w.r.t.\ the Lebesgue measure is bounded above and below:
\[
\exists\ 0 < m_1 \leq m_2 < \infty \quad \forall x \in B \qquad m_1 \lambda (\dx{x}) \leq \dx{\mu} (x) \leq m_2 \lambda (\dx{x})
.\]
Let $\phi \colon \RR^n \to \RR_+$ be a function such that:
\begin{enumerate}[(i)]
\item $\phi$ is zero outside $B$;
\item There exist $0 < c_1 \leq c_2 < \infty$ such that for all $x \in B$, $c_1 \dist (x, \C{B})^2 \leq \phi(x) \leq c_2 \dist(x, \C{B})^2$.
\item $\phi$ is $k$-Lipschitz for some $k < \infty$.
\end{enumerate}
Then, denoting $a \coloneqq \norm{\phi \cdot \nu}_1 \div \norm{\phi \cdot \mu}_1$,
\[\label{eqn:localisation}
W_2 (a \phi \cdot \mu, \phi \cdot \nu) \leq C (n)^{1/2} \frac{c_2^{3/2} m_2^{3/2}}{c_1^{3/2} m_1^{3/2}} k c_1^{-1/2} W_2 (\mu, \nu)
,\]
for $C (n) < \infty$ some absolute constant only depending on $n$. Moreover, taking $C (n) \coloneqq 2^{11} n$ fits.%
\footnote{Though of course the factor $2^{11}$ may be strongly suboptimal.}
\end{theorem}

\begin{remark}
Actually the constraint that the support of $\phi$ is a ball is of little importance: we could assume as well that it would be a cube, a simplex, or many other shapes, as the corollary below shows:

\begin{corollary}\label{cor:localisation.othershapes}
Make the same assumptions as in Theorem~\ref{thm:localisation.sphere}, except that $B$ need not be a ball: instead, we only assume that, denoting by $B_{\circ}$ the (true) ball having the same volume as $B$, there exists a bijection $\Phi \colon B \leftrightarrow B_{\circ}$ mapping the uniform measure on $B$ onto the uniform measure on $B_{\circ}$ (i.e.\ such that $\Phi \pushfwd (\1{B} \cdot \lambda) = \1{B_{\circ}} \cdot \lambda$) such that $\Phi$ is bi-Lipschitz (i.e.\ such that both $\Phi$ and $\Phi^{-1}$ are Lipschitz). Denote by $\norm{\Phi}_{\Lipschitz}$ and $\norm{\Phi^{-1}}_{\Lipschitz}$ the optimal Lipschitz constants for resp.\ $\Phi$ and $\Phi^{-1}$. Then, the conclusion of Theorem~\ref{thm:localisation.sphere} remains true, except that now you have to replace the factor $C (n)$ by
\[
(\norm{\Phi}_{\Lipschitz} \norm{\Phi^{-1}}_{\Lipschitz})^{10} C (n)
.\]
\end{corollary}

\begin{proof}
Consider the measures $\mu_{\circ} \coloneqq \Phi \pushfwd \mu$ and $\nu_{\circ} \coloneqq \Phi \pushfwd \nu$, and the bump function $\phi_{\circ} \coloneqq \phi \circ \Phi^{-1}$; then, $\mu_{\circ}$, $\nu_{\circ}$ and $\phi_{\circ}$ satisfy the original assumptions of Theorem~\ref{thm:localisation.sphere}, the roles of \squo{$m_1$} and \squo{$m_2$} (in the ball situation) being held by $m_1$ and $m_2$ (in the general situation) themselves, the role of \squo{$k$} being held by $\norm{\Phi^{-1}}_{\Lipschitz} k$, and the roles of \squo{$c_1$} and \squo{$c_2$} being held by $c_1 \div \norm{\Phi}_{\Lipschitz}^2$ and $c_2 \norm{\Phi^{-1}}_{\Lipschitz}^2$. Therefore, applying \eqref{eqn:localisation}:
\[
W_2 (a \phi_{\circ} \cdot \mu_{\circ}, \phi_{\circ} \cdot \nu_{\circ}) \leq C (n)^{1/2} \norm{\Phi}_{\Lipschitz}^4 \norm{\Phi^{-1}}_{\Lipschitz}^4 \frac{c_2^{3/2} m_2^{3/2}}{c_1^{3/2} m_1^{3/2}} W_2 (a \mu_{\circ}, \nu_{\circ})
.\]
But the optimal transportation plan from $a \mu$ to $\nu$, with cost $W_2 (\mu, \nu)^2$, can be pushed forward by $\Phi$ into a (not optimal in general) transportation plan from $a \mu_{\circ}$ to $\nu_{\circ}$, whose cost will then be $\leq \norm{\Phi}_{\Lipschitz}^2 W_2 (\mu, \nu)^2$; so $W_2 (a \mu_{\circ}, \nu_{\circ}) \leq \norm{\Phi}_{\Lipschitz} W_2 (a \mu, \nu)$. Similarly $W_2 (a \phi \cdot \mu, \phi \cdot \nu) \leq \norm{\Phi^{-1}}_{\Lipschitz} W_2 (a \phi_{\circ} \cdot \mu_{\circ}, \phi_{\circ} \cdot \nu_{\circ})$. The announced result follows.
\end{proof}
\end{remark}

\subsection{Proof of the main theorem}

In the sequel we will shorthand $W_2 (\mu, \nu) \eqqcolon w$, and also $\phi \cdot \mu \eqqcolon \hat{\mu}$, resp.\ $\phi \cdot \nu \eqqcolon \hat{\nu}$. Let $g \eqqcolon \id + S$ be a map achieving optimal transportation from $\nu$ to $\mu$, i.e.\ such that $\mu = g \pushfwd \nu$ with $\int_{\RR^n} \abs{S (y)}^2 \dx{\nu} (y) = w^2$.%
\footnote{\label{ftn:f}Actually such an $g$ does not always exist, as it can occur that the optimal transportation plan from $\nu$ to $\mu$ \dquo{splits points} if $\nu$ is not regular enough. However it would suffice to use the general formalism of transportation plans to handle that case: we do not do it here to keep notation light, but this is straightforward. Also note that it is not obvious that the infimum in \eqref{eqn:W2} is attained: again, that is not a real problem as our proof still works by considering a sequence of transportation plans approaching optimality.}

Our strategy will consist in transforming $\hat{\nu}$ into $a \hat{\mu}$ according to the following procedure:
\begin{enumerate}
\item[\stepone]
We apply the transportation plan $g$ to $\hat{\nu}$; this transforms $\hat{\nu}$ into some measure $\hat{\mu}^*$. The measure $\hat{\mu}^*$ is not supported by $B$ \emph{a priori}, so we split it into $\hat{\mu}^*_B + \hat{\mu}^*_{\csf} \coloneqq \1{B} \cdot \hat{\mu}^* + \1{\C{B}} \cdot \hat{\mu}^*$.
\item[\steptwo]
Denoting $a_{\csf} \coloneqq \norm{\hat{\mu}^*_{\csf}}_1 \div \norm{\hat{\mu}}_1$, we then transform $\hat{\mu}^*_{\csf}$ into $a_{\csf} \hat{\mu}$ according to an arbitrary transference plan;
\item[\stepthree]
Finally, denoting $a_B \coloneqq \norm{\hat{\mu}^*_B}_1 \div \norm{\hat{\mu}}_1$,%
\footnote{Observe that $a_B + a_{\csf} = a$.}
we transform $\hat{\mu}^*_B$ into $a_B \hat{\mu}$ according to the optimal transference plan: the cost of this operation is $W_2 (\hat{\mu}^*_B, a_B \hat{\mu})$, which we bound above by $2 \norm{\hat{\mu}^*_B - a_B \hat{\mu}}_{\dHmo (a_B \hat{\mu})}$ thanks to Theorem~\ref{thm:W2<dHmo}.
\end{enumerate}
Then, denoting by $W_2 (\stepone), W_2 (\steptwo), W_2 (\stepthree)$ the respective Wasserstein distances of these steps, we shall have $W_2 (\hat{\nu}, a\hat{\mu}) \leq W_2 (\stepone) + (W_2 (\steptwo)^2 + W_2(\stepthree)^2)^{1/2}$.

\bigskip

Let us begin with bounding the cost of Step~\stepone. The squared cost of this step is
\begin{multline}
W_2 (\stepone)^2 = \int \abs{S (y)}^2 \dx{\hat{\nu}} (y) = \int \abs{S (y)}^2 \phi(y) \dx{\nu} (y) \\
\leq \sup \phi \times \int \abs{S (y)}^2 \dx{\nu} (y) = \sup \phi \times w^2 \leq c_2 R^2 w^2
,\end{multline}
whence $W_2(\stepone) \leq c_2^{1/2} R w$.

\medskip

Now consider Step~\steptwo. As $a_{\csf} \hat{\mu}$ is supported by $B$, one has obviously
\[
W_2 (\steptwo)^2 \leq \int \paren[\big]{\dist (x, B) + 2 R}^2 \dx{\hat{\mu}^*_{\csf}} (x) = \int_{\C{B}} \paren[\big]{\dist (x, B) + 2 R}^2 \dx{\hat{\mu}^*} (x)
.\]
From that we deduce that $W_2(\steptwo) \leq 2 c_2^{1/2} R w$ by the following computation:
\begin{multline}
\int_{\C{B}} \paren[\big]{\dist (x, B) + 2 R}^2 \dx{\hat{\mu}^*} (x)
= \int_{g (y) \notin B} \paren[\big]{\dist (g (y), B) + 2 R}^2 \phi(y) \dx{\nu} (y) \\
\leq c_2 \int_{\substack{y \in B\\ g (y) \notin B}} \paren[\big]{\dist (g (y), B) + 2 R}^2 \dist (y, \C{B})^2 \dx{\nu} (y) \\
\leq c_2 \int_{\substack{y \in B\\ g (y) \notin B}} \paren[\big]{R \dist (g (y), B) + 2 R \dist (y, \C{B})}^2 \dx{\nu} (y) \\
\leq 4 c_2 R^2 \int_{\substack{y \in B\\ g (y) \notin B}} \paren[\big]{\dist (g (y), B) + \dist (y, \C{B})}^2 \dx{\nu} (y) \\
\leq 4 c_2 R^2 \int \abs{y - g (y)}^2 \dx{\nu} (y)
= 4 c_2 R^2 w^2
.\end{multline}

Step~\stepthree\ is the difficult one. We begin with observing that it is easy to bound the $L^2 (B)$ distance between $\hat{\mu}^*_B$ and $\hat{\mu}$: indeed, denoting by $f \eqqcolon \id + T$ the inverse map of $g$%
\footnote{For $f$ to exist, $g$ should be bijective, which is not always true \emph{stricto sensu}; but we can safely carry out the reasoning with pretending so, by the same argument as in Footnote~\ref{ftn:f} on page~\pageref{ftn:f}.},
\begin{multline}\label{easy-bound_L2muB}
\norm{\hat{\mu}^*_B - \hat{\mu}}_{L^2 (\1{B} \cdot \mu)}^2
= \int_B \paren[\bigg]{\frac{\dx{\hat{\mu}^*} (x) - \phi (x) \dx{\mu} (x)}{\dx{\mu} (x)}}^2 \dx{\mu} (x) \\
= \int_B \paren[\big]{\phi (f (x)) - \phi (x)}^2 \dx{\mu} (x) \\
\leq k^2 \int_{\RR^n} \abs{x - f (x)}^2 \dx{\mu} (x) = k^2 \int \abs{T (x)}^2 \dx{\mu} (x) = k^2 w^2
,\end{multline}
(where we used that $\dx{\hat{\mu}^*} (x) = \dx{\hat{\nu}} (f (x)) = \phi (f (x)) \dx{\nu} (f (x)) = \phi (f (x)) \dx{\mu} (x)$), so that
\[\label{easybound_L2B}
\norm{\hat{\mu}^*_B - \hat{\mu}}_{L^2 (B)}^2 \leq k^2 m_2 w^2
.\]

Now we have to link $\norm{\placeholder}_{L^2 (B)}$ with $\norm{\placeholder}_{\dHmo (\mu)}$. This is achieved by the following lemma, whose proof is postponed:
\begin{lemma}\label{lem:dHmo<L2B}
Define $\hat{\lambda}$ to be the measure on $B$ such that $\hat{\lambda} (\dx{x}) \coloneqq \dist(x, \C{B})^2 \lambda (\dx{x})$. Then, for any signed measure $m$ on $B$ having total mass zero:
\[
\norm{m}_{\dHmo (\hat{\lambda})} \leq C_1 (n)^{1/2} \norm{m}_{L^2 (B)}
,\]
where $C_1 (n)$ is some absolute constant only depending on $n$. Moreover, taking $C_1 (n) \coloneqq \paren[\big]{(2 e + 1) n - 1} \vee 8 e$ fits.
\end{lemma}

Thanks to Theorem~\ref{thm:W2<dHmo} and Lemma~\ref{lem:dHmo<L2B}, we have that
\begin{multline}\label{eqn:stepthree<L2B}
W_2 (\stepthree) \leq 2 \norm{a_B \hat{\mu} - \hat{\mu}^*_B}_{\dHmo (a_B \hat{\mu})}
\leq 2 (a_B c_1 m_1)^{-1/2} \norm{a_B \hat{\mu} - \hat{\mu}^*_B}_{\dHmo (\hat{\lambda})} \\
\leq 2 C_1 (n)^{1/2} (a_B c_1 m_1)^{-1/2} \norm{a_B \hat{\mu} - \hat{\mu}^*_B}_{L^2 (B)}
.\end{multline}
Next, we compute
\begin{multline}\label{eqn:L2B<W2}
\norm{a_B \hat{\mu} - \hat{\mu}^*_B}_{L^2 (B)}
= \norm[\Big]{\tsfrac{\norm{\hat{\mu}^*_B}_1}{\norm{\hat{\mu}}_1} \hat{\mu} - \hat{\mu}^*_B}_{L^2 (B)}
\leq \tsfrac{\abs{\norm{\hat{\mu}^*_B}_1 - \norm{\hat{\mu}}_1}}{\norm{\hat{\mu}}_1} \norm{\hat{\mu}}_{L^2 (B)} + \norm{\hat{\mu}^*_B - \hat{\mu}}_{L^2 (B)} \\
\leq \tsfrac{\norm{\hat{\mu}}_{L^2 (B)}}{\norm{\hat{\mu}}_1} \norm{\hat{\mu}^*_B - \hat{\mu}}_1 + \norm{\hat{\mu}^*_B - \hat{\mu}}_{L^2 (B)}
\leq \paren[\Big]{\tsfrac{\norm{\hat{\mu}}_{L^2 (B)}}{\norm{\hat{\mu}}_1} \lambda (B)^{1/2} + 1} \norm{\hat{\mu}^*_B - \hat{\mu}}_{L^2 (B)} \\
\leq \paren[\Big]{\tsfrac{c_2 m_2}{c_1 m_1} \tsfrac{\lambda (B)^{1/2} \norm{\hat{\lambda}}_{L^2 (B)}}{\norm{\hat{\lambda}}_1} + 1} \norm{\hat{\mu}^*_B - \hat{\mu}}_{L^2 (B)}
\bstackrel{\footnotemark}{\leq} \sqrt{6} \tsfrac{c_2 m_2}{c_1 m_1} \norm{\hat{\mu}^*_B - \hat{\mu}}_{L^2 (B)} \\
\bstackrel{\eqref{easybound_L2B}}{\leq} \sqrt{6} \tsfrac{c_2 m_2}{c_1 m_1} k m_2^{1/2} w
,\end{multline}
\footnotetext{This step comes from the computation $\lambda (B)^{1/2} \norm{\hat{\lambda}}_{L^2 (B)} \div \norm{\hat{\lambda}}_1 = {(\int_0^1 r^{n - 1} \dx{r})^{1/2}} \* {\paren[\big]{\int_0^1 (1 - r)^4 r^{n - 1} \dx{r}}^{1/2}} \div {\paren[\big]{\int_0^1 (1 - r)^2 r^{n - 1} \dx{r}}} = \paren[\big]{6 (1 + n) (2 + n) \div (3 + n) (4 + n)}^{1/2} \leq \sqrt{6}\ \forall n$.}
so that, combining~\eqref{eqn:stepthree<L2B} and~\eqref{eqn:L2B<W2}, we have got:
\[\label{eqn:stepthree.raw}
W_2 (\stepthree) \leq 2 \sqrt{6} C_1 (n)^{1/2} a_B^{-1/2} \frac{c_2 m_2^{3/2}}{c_1 m_1^{3/2}} \frac{k}{c_1^{1/2}} w
.\]

Equation~\eqref{eqn:stepthree.raw} is the kind of bound we were looking for, \emph{provided $a_B \lesssim 1$}. Though this will be the case in practice (since we are mainly interested in cases where $\nu$ is close to $\mu$ and thus $\hat{\mu}^*$ is close to $\hat{\mu}$), this is not quite satisfactory yet. So, what can we do when $a_B \ll 1$, that is, when $\norm{\hat{\mu}^*_B}_1 \ll \norm{\hat{\mu}}_1$? In fact that case is easier, because transportation between small measures has low cost, while $w$ has to be large to make $\hat{\mu}^*_B$ very different from $\hat{\mu}$.

The computations are the following. First, it is obvious that
\[\label{eqn:stepthree.trivial}
W_2 (\stepthree) = W_2 (\hat{\mu}^*_B, a_B \hat{\mu}) \leq 2 R \norm{\hat{\mu}^*_B}_1^{1/2}
.\]
Next, observing that $\phi (f (x)) \geq \frac{c_1}{c_2} \phi (x) - 2 c_1 \dist (x, \C{B}) \abs{T (x)}$,%
\footnote{This follows from the computation:
\begin{multline}
\phi (f (x)) \geq c_1 \dist (f (x), \C{B})^2 \geq c_1 \paren[\big]{\dist (x, \C{B}) - \abs{T (x)}}_+^2 \\
\geq c_1 \dist (x, \C{B})^2 - 2 c_1 \dist (x, \C{B}) \abs{T (x)} \geq \frac{c_1}{c_2} \phi (x) - 2 c_1 \dist (x, \C{B}) \abs{T (x)}
.\end{multline}}
we compute that
\begin{multline}\label{eqn:W2>tildemu-xiB}
\norm{\hat{\mu}^*_B}_1
= \int_B \phi (f (x)) \dx{\mu} (x)
\geq \int_B \paren[\Big]{\tsfrac{c_1}{c_2} \phi (x) - 2 c_1 \dist (x, \C{B}) \abs{T (x)}} \dx{\mu} (x) \\
\geq \tsfrac{c_1}{c_2} \norm{\hat{\mu}}_1 - 2 c_1 \paren[\Big]{\int_B \dist (x, \C{B})^2 \dx{\mu} (x)}^{1/2} \paren[\Big]{\int_B \abs{T (x)}^2 \dx{\mu} (x)}^{1/2} \\
= \tsfrac{c_1}{c_2} \norm{\hat{\mu}}_1 - 2 c_1 \norm{\dist (\placeholder, \C{B})^2 \cdot \mu}_1^{1/2} w
\geq \tsfrac{c_1}{c_2} \norm{\hat{\mu}}_1 - 2 c_1 m_2^{1/2} \norm{\hat{\lambda}}_1^{1/2} w
,\end{multline}
whence
\[
w \geq \frac{\paren[\big]{\tsfrac{c_1}{c_2} \norm{\hat{\mu}}_1 - \norm{\hat{\mu}^*_B}_1}_+}{2 c_1 m_2^{1/2} \norm{\hat{\lambda}}_1^{1/2}} = \frac{\paren[\big]{\tsfrac{c_1}{c_2} - a_B}_+ \norm{\hat{\mu}}_1}{2 c_1 m_2^{1/2} \norm{\hat{\lambda}}_1^{1/2}} \geq \frac{m_1^{1/2}}{2 c_1 m_2^{1/2}} \paren[\big]{\tsfrac{c_1}{c_2} - a_B}_+ \norm{\hat{\mu}}_1^{1/2}
.\]
So,
\[\label{eqn:stepthree.variant}
W_2 (\stepthree)
\leq 2 R \norm{\hat{\mu}^*_B}_1^{1/2}
= 2 R a_B^{1/2} \norm{\hat{\mu}}_1^{1/2}
\leq 4 R c_1^{1/2} \tsfrac{m_2^{1/2}}{m_1^{1/2}} \frac{a_B^{1/2}}{(\frac{c_1}{c_2} - a_B)_+} w
.\]

In the end, choosing either \eqref{eqn:stepthree.raw} if $a_B \geq c_1 \div 2 c_2$ or \eqref{eqn:stepthree.variant} if $c_1 \div 2 c_2$, and observing that $c_1 \leq k R^{-1}$, one has always:
\[\label{bound:W2(3)}
W_2 (\stepthree) \leq (4 \sqrt{3} C_1 (n)^{1/2} \vee 4 \sqrt{2}) \frac{c_2^{3/2} m_2^{3/2}}{c_1^{3/2} m_1^{3/2}} \frac{k}{c_1^{1/2}} w
.\qed\]

\begin{remark}
To bound $W_2 (\stepthree)$ in the situation where $a_B \ll 1$, we could also have started from \quotes{$\phi (f (x)) \geq \phi (x) - k \abs{T (x)}$} (instead of \quotes{$\phi (f (x)) \geq \frac{c_1}{c_2} \phi (x) - 2 c_1 \dist (x, \C{B}) \abs{T (x)}$}) to get another bound analogous to \eqref{eqn:W2>tildemu-xiB}. Following such an approach, the factor $(c_2 \div c_1)^{3/2}$ in \eqref{eqn:stepthree.variant} would be improved into $(c_2 \div c_1)$ in the analogous formula; however the dimensional factor would behave in $O (n)$ rather than in $O (n^{1/2})$.
\end{remark}

\subsection{Proof of Lemma~\ref{lem:dHmo<L2B}}

It still remains to prove Lemma~\ref{lem:dHmo<L2B}, whose statement we recall to be:
\begin{lemma*}
Denoting $\hat{\lambda} \coloneqq \dist(\placeholder, \C{B})^2 \cdot \lambda$, one has, for any signed measure $m$ on $B$ having total mass zero:
\[\label{dHmo<L2B.explicite}
\norm{m}_{\dHmo (\hat{\lambda})} \leq \paren[\big]{\paren[\big]{(2 e + 1) n - 1} \vee 8 e}^{1/2} \norm{m}_{L^2 (B)}
.\]
—In the sequel, \quotes{$\paren[\big]{(2 e + 1) n - 1} \vee 8 e$} will be shorthanded into \quotes{$C_1 (n)$}.
\end{lemma*}

\begin{remark}
The bound \eqref{dHmo<L2B.explicite} is within a constant factor of being optimal, uniformly in $n$, as one sees by taking $f (x) = x_z$ in \eqref{eqn:6921} ($x_z$ denoting the $z$-coordinate of $x$).
\end{remark}

\begin{proof}[\proofname\ of the lemma]
We begin with translating the lemma into a functional analysis statement by a duality argument. Recall the duality definition of $\norm{m}_{\dHmo (\hat{\lambda})}$ from \S~\ref{sec:equivalence}:
\[\label{eqn:dHmo.dual}
\norm{m}_{\dHmo (\hat{\lambda})} \coloneqq \sup \setst{\abs{\inner{f}{m}}}{\norm{f}_{\dot{H}^1 (\hat{\lambda})} \leq 1}
.\]
There is a similar duality formula for~$\norm{m}_{L^2 (B)}$:
\[\label{eqn:L2B.dual}
\norm{m}_{L^2 (B)} = \sup \setst{\abs{\inner{f}{m}}}{\norm{f}_{L^2 (B)} \leq 1}
,\]
where, for $f$ a \emph{function}, $\norm{f}_{L^2 (B)}$ has its usual meaning, namely $\norm{f}_{L^2 (B)} \coloneqq \paren[\big]{\int_B f (x)^2 \dx{\lambda} (x)}^{1/2}$.
Since $m$ is assumed to have total mass zero, $\abs{\inner{f}{m}}$ does not change when one adds a constant to $f$. On the other hand, when $f$ describes the set $\setst{\norm{f_0 + a}}{a \in \RR}$, $\norm{f}_{L^2 (B)}$ is minimal when $a$ is such that $f$ has zero mean on $B$, while the value of $\norm{f}_{\dot{H}^{1} (\hat{\lambda})}$ remains constant.%
\footnote{Here we implicitly assume that $\int_B \abs{f (x)} \dx{\lambda} (x)$, which is legit since an approximation argument allows to restrict the suprema in \eqref{eqn:dHmo.dual} and \eqref{eqn:L2B.dual} to those $f$ having a $\mcal{C}^{\infty}$ continuation on $\closure (B)$.}
As a consequence, we can restrict the supremum in \eqref{eqn:dHmo.dual} and \eqref{eqn:L2B.dual} to those $f$ \emph{having zero mean on $B$}. Thus, the lemma will be implied%
\footnote{Actually there is even equivalence.}
by proving that
\[\label{eqn:6921}
\inner{f}{\1{B} \cdot \lambda} = 0 \quad \Rightarrow \quad \norm{f}_{L^2 (B)} \leq C_1 (n)^{1/2} \norm{f}_{\dot{H}^{1} (\hat{\lambda})}
.\]
Going back to the definitions of $\norm{\placeholder}_{\dHmo (\hat{\lambda})}$ and $\norm{\placeholder}_{L^2 (B)}$, relaxing the condition on $f$ to be centred by projecting it orthogonally in $L^2 (B)$ onto the subspace of centred functions, and denoting by $P$ the uniform probability measure on $B$, Equation~\eqref{eqn:6921} turns into:
\[\label{eqn:Poincare}
\forall f \qquad \Var_{P} (f) \leq C_1 (n) \int \dist (x, \C{B})^2 \abs{\nabla{f} (x)}^2 \dx{P} (x)
,\]
which we recognize to be a weighted Poincaré inequality.

To prove \eqref{eqn:Poincare}, the first key idea (inspired by \citep{Bobkov}) is to separate radial and spherical coordinates. This is, considering the bijection
\[
\begin{aligned}
\phi \colon (0, R) \times \sphere{n - 1} & \to B \setminus \{0\} \\
(r, \theta) & \mapsto r \theta
\end{aligned}
\]
(the origin of space being set at the center of $B$), we introduce the measure $\tilde{P} \coloneqq \phi^{-1} \pushfwd P$, which is obviously the product measure $\tilde{P}_r \otimes \tilde{P}_{\theta}$, where $\tilde{P}_r$ is the probability measure on $(0, R)$ such that $\dx{\tilde{P}_r} (r) \coloneqq n R^{-n} r^{n - 1} \dx{r}$, resp.\ $\tilde{P}_{\theta}$ is the uniform measure on the sphere $\sphere{n - 1}$.
With this notation, we perform can a change of variables to see that \eqref{eqn:Poincare} is equivalent to proving that, for all $g \in L^2 (\tilde{P})$:
\[\label{eqn:Poincare.rtheta}
C_1 (n)^{-1} \Var_{\tilde{P}} (g) \leq \int_0^R \int_{\sphere{n - 1}} (R - r)^2 \bigl(\abs{\nabla_r{g} (r, \theta)}^2 + r^{-2} \abs{\nabla_{\theta}{g} (r, \theta)}^2\bigr) \dx{\tilde{P}_r} (r) \dx{\tilde{P}_{\theta}} (\theta)
,\]
where $\nabla_r$ and $\nabla_{\theta}$ denote the gradient along resp.\ the $r$ coordinate and the $\theta$ coordinate.%
\footnote{In the latter case, we have to use the Riemannian definition of the gradient on $\sphere{n - 1}$.}
We will denote the right-hand side of \eqref{eqn:Poincare.rtheta} by $\mcal{E} (g, g)$.

Because $\tilde{P} = \tilde{P}_r \otimes \tilde{P}_{\theta}$, we know that $L^2 (\tilde{P})$ can be seen as (the closure of) the tensor product of $L^2 (\tilde{P}_r)$ and $L^2 (\tilde{P}_{\theta})$:
\[\label{eqn:L2tP=tensprod}
L^2 (\tilde{P}) = \closure (L^2 (\tilde{P}_r) \overset{\perp}{\otimes} L^2 (\tilde{P}_{\theta}))
,\]
where the symbol \squo{$\overset{\perp}{\otimes}$} means that the Hilbertian structure of $L^2 (\tilde{P})$ is compatible with the Hilbertian structures of $L^2 (\tilde{P}_r)$ and $L^2 (\tilde{P}_{\theta})$—i.e., that $\inner{h_{\mrm{a}} \otimes u_{\mrm{a}}}{h_{\mrm{b}} \otimes u_{\mrm{b}}}_{L^2 (\tilde{P})} = \inner{h_{\mrm{a}}}{h_{\mrm{b}}}_{L^2 (\tilde{P}_r)} \times \inner{u_{\mrm{a}}}{u_{\mrm{b}}}_{L^2 (\tilde{P}_{\theta})}$.
Now consider the spherical harmonics $Y_0, Y_1, \ldots$, which by definition are an orthonormal basis, in $L^2 (\tilde{P}_{\theta})$, of eigenfunctions of the Laplace–Beltrami operator $\Delta$ on $\sphere{n - 1}$; and call $\ell_0, \ell_1, \ldots$ the associated eigenvalues, which are known to be such that (up to permuting indices) $Y_0 \equiv 1$ with $\ell_0 = 0$, and $\ell_i \leq -(n - 1)\ \forall i \neq 0$ (see for instance \citep{sphericalharmonics}).
By construction, $L^2 (\tilde{P}_{\theta}) = \closure \paren[\big]{\mathop{\overset{\perp}{\bigoplus}}_{i \in \NN} (\RR \cdot Y_i)}$; therefore, one has that
\[
L^2 (\tilde{P}) = \closure \paren[\Big]{\mathop{\overset{\perp}{\bigoplus}}_{i \in \NN} L^2 (\tilde{P}_r) \cdot Y_i}
:\]
in other words, the functions of $L^2 (\tilde{P})$ are those of the form
\[
g (r, \theta) = \sum_{i \in \NN} h_i (r) Y_i (\theta)
,\]
with $\sum_i \norm{h_i}_{L^2 (\tilde{P}_r)}^2 < \infty$, and the correspondence is bijective. An interesting point is that, then, one has:
\[
\Var_{\tilde{P}} (g) = \Var_{\tilde{P}_r} (h_0) + \sum_{i \neq 0} \norm{h_i}_{L^2 (\tilde{P}_r)}^2
.\]

On the other hand, one has
\[
\mcal{E} (g, g) = -\langle L g, g \rangle_{L^2 (\tilde{P})}
,\]
where
\[\label{eqn:L}
\paren[\big]{L g} (r, \theta) \coloneqq (R - r)^2 \Delta_r g + \paren[\bigg]{(n - 1) \frac{(R - r)^2}{r} - 2 (R - r)} \mathbf{e}_r \cdot \nabla_r{g} + \frac{(R - r)^2}{r^2} \Delta_{\theta} g
.\]
From \eqref{eqn:L} we see that, since the $Y_i$ are eigenfunctions of $\Delta_\theta$, all the $L^2 (\tilde{P}_r) \cdot Y_i$ are invariant by $L$, and that one has:
\[\label{eqn:Egg.alphai}
\mcal{E} (g, g) = \sum_{i \in \NN} \int_0^R \paren[\bigg]{(R - r)^2 \abs{\nabla{h_i} (r)}^2 - \ell_i \frac{(R - r)^2}{r^2} h_i (r)^2} \tilde{P}_r (\dx{r})
.\]

So, proving \eqref{eqn:Poincare.rtheta} becomes equivalent to proving that both following formulas hold for all $h \in L^2 (\tilde{P}_r)$:
\[\label{eqn:Poincare.i=0}
\Var_{\tilde{P}_r} (h) \leq C_1 (n) \int_0^R (R - r)^2 \abs{\nabla{h} (r)}^2 \tilde{P}_r (\dx{r})
;\]
\[\label{eqn:Poincare.ineq0}
\norm{h}_{L^2 (\tilde{P}_r)}^2 \leq C_1 (n) \int_0^R \paren[\bigg]{(R - r)^2 \abs{\nabla{h} (r)}^2 + (n - 1) \frac{(R - r)^2}{r^2} h (r)^2} \tilde{P}_r (\dx{r})
.\]

\bigskip

Let us start with \eqref{eqn:Poincare.i=0}. In all the sequel of the proof, we introduce the following notation:
\begin{align}
w (r) & \coloneqq (R - r)^{-3/2}; \\
b & \coloneqq 1 - n^{-1}.
\end{align}

By the Cauchy–Schwarz inequality, one has, for all $r \in (b R, R)$:
\[
\paren[\big]{h (r) - h (b R)}^2 = \paren[\Big]{\int_{b R}^r h' (s) \dx{s}}^2 \leq \paren[\Big]{\int_{b R}^r w (s) \dx{s}} \times \int_{b R}^r w (s)^{-1} \abs{\nabla{h} (s)}^2 \dx{s}
,\]
which yields, with the chosen form of $w (\placeholder)$:
\begin{multline}\label{CSbound.concrete}
(h (r) - h (b R))^2 \leq 2 \paren[\big]{(R - r)^{-1/2} - (R - b R)^{-1/2}} \int_{b R}^r (R - s)^{3/2} \abs{\nabla{h} (s)}^2 \dx{s} \\
\leq 2 (R - r)^{-1/2} \int_{b R}^r (R - s)^{3/2} \abs{\nabla{h} (s)}^2 \dx{s}
.\end{multline}
Integrating and using Fubini's formula, it follows that
\begin{multline}\label{eqn:Poincare-boundary}
\int_{b R}^R \paren[\big]{h (r) - h (b R)}^2 \dx{\tilde{P}_r} (r) \leq \\
2 \int_{s = b R}^R \Bigl(\int_{r = s}^R n R^{-n} (R - r)^{-1/2} r^{n - 1} \dx{r}\Bigr) (R - s)^{3/2} \abs{\nabla{h} (s)}^2 \dx{s} \\
\leq 2 \int_{s = b R}^R \Bigl(\int_{r = s}^R n R^{-n} (b^{-1} s)^{n - 1} (R - r)^{-1/2} \dx{r}\Bigr) (R - s)^{3/2} \abs{\nabla{h} (s)}^2 \dx{s} \\
= 2 b^{-(n - 1)} \int_{s = b R}^R \Bigl(\int_{r = s}^R (R - r)^{-1/2} \dx{r}\Bigr) (R - s)^{3/2} \abs{\nabla{h} (s)}^2 \dx{\tilde{P}_r} (s) \\
= 4 b^{-(n - 1)} \int_{s = b R}^R (R - s)^2 \abs{\nabla{h} (s)}^2 \dx{s}.
\end{multline}

One can apply the same line of reasoning for $r \in (0, b R)$: the (unweighted this time) Cauchy–Schwarz inequality then yields $\paren[\big]{h (r) - h (b R)}^2 \leq (b R - r) \* \int_r^{b R} \abs{\nabla h (s)}^2 \dx{s}$, whence:
\begin{multline}\label{eqn:Poincare-centre}
\int_0^{b R} \paren[\big]{h (r) - h (b R)}^2 \dx{\tilde{P}_r} (r) \leq
\int_{s = 0}^{b R} \paren[\Big]{\int_{r = 0}^s n R^{-n} (b R - r) r^{n - 1} \dx{r}} \abs{\nabla{h} (s)}^2 \dx{s} \\
\leq R^{- (n - 1)} \int_{s = 0}^{b R} \paren[\Big]{\int_{r = 0}^s n r^{n - 1} \dx{r}} \abs{\nabla{h} (s)}^2 \dx{s}
= R \int_0^{b R} \abs{\nabla{h} (s)}^2 s^n \dx{s} \\
\leq n^{-1} R^2 \int_0^{b R} \abs{\nabla{h} (s)}^2 \dx{\tilde{P}_r} (s)
\leq n^{-1} (1 - b)^{-2} \int_0^{b R} (R - s)^2 \abs{\nabla{h} (s)}^2 \dx{\tilde{P}_r} (s)
.\end{multline}

Summing \eqref{eqn:Poincare-boundary} and \eqref{eqn:Poincare-centre}, we get that
\[\label{eqn:spectral-gap}
\int_0^R \paren[\big]{h (r) - h (b R)}^2 \dx{\tilde{P}_r} (r) \leq \paren[\big]{4 b^{-(n - 1)} \vee n^{-1} (1 - b)^{-2}} \int_0^s (R - s)^2 \abs{\nabla{h} (s)}^2 \dx{\tilde{P}_r} (s)
,\]
where $\paren[\big]{4 b^{-(n - 1)} \vee n^{-1} (1 - b)^{-2}}$ can itself be bounded by $(n \vee 4e)$. The left-hand-side of \eqref{eqn:spectral-gap} being an upper bound for $\Var_{\tilde{P}_r} (h)$, this proves \eqref{eqn:Poincare.i=0}.

\bigskip

Now we turn to \eqref{eqn:Poincare.ineq0}. For $r \in (b R, R)$ we have, similarly to \eqref{CSbound.concrete}, that
\[
\paren[\big]{h (r) - h (b r)}^2 \leq 2 (R - r)^{-1/2} \int_{b r}^r (R - s)^{3/2} \abs{\nabla h (s)}^2 \dx{s}
,\]
so that
\[
h (r)^2 \leq 2 h (b r)^2 + 4 (R - r)^{-1/2} \int_{b r}^r (R - s)^{3/2} \abs{\nabla h (s)}^2 \dx{s}
.\]
Then, integrating and applying Fubini's formula:
\begin{multline}\label{eqn:3954}
\int_{b R}^R h (r)^2 \dx{\tilde{P}_r} (r) \leq 2 \int_{b R}^R h (b r)^2 \dx{\tilde{P}_r} (r) + \null \\ 4 \int_{s = b^2 R}^R \paren[\Big]{\int_{r = s \vee b R}^{b^{-1} s \wedge R} n R^{-n} r^{n - 1} (R - r)^{-1/2} \dx{r}} (R - s)^{3/2} \abs{\nabla h (s)}^2 \dx{s}
.\end{multline}

By change of variables, the first term of the right-hand side of \eqref{eqn:3954} is equal to $2 b^{-n} \int_{b^2 R}^{b R} h (s)^2 \dx{\tilde{P}_r} (s)$, which we can bound by
\begin{multline}\label{bound-lambda0.0}
2 b^{-(n - 2)} \frac{(1 - b)^{-2}}{n - 1} \int_{b^2 R}^{b R} (n - 1) \frac{(R - r)^2}{r^2} h (s)^2 \dx{\tilde{P}_r} (s) \\
\leq 2 n e \int_0^R (n - 1) \frac{(R - r)^2}{r^2} h (s)^2 \dx{\tilde{P}_r} (s)
.\end{multline}
The second term of the right-hand side of \eqref{eqn:3954} is itself bounded by
\begin{multline}\label{bound-lambda0.1}
4 b^{-(n - 1)} \int_{s = b^2 R}^R \paren[\Big]{\int_{r = s}^R (R - r)^{-1/2} \dx{r}} (R - s)^{3/2} \abs{\nabla h (s)}^2 \dx{\tilde{P}_r} (s) \\
\leq 8 e \int_0^R (R - s)^2 \abs{\nabla h (s)}^2 \dx{\tilde{P}_r} (s)
.\end{multline}
This way, we have bounded $\int_{b R}^R h (r)^2 \dx{\tilde{P}_r} (r)$.

On the other hand, it is trivial that, for $r \leq b R$,
\[
h (r)^2 \leq \frac{b^2}{(n - 1) (1 - b)^2} \times (n - 1) \frac{(R - r)^2}{r^2} h (r)^2
,\]
whence:
\[\label{bound-lambda0.2}
\int_0^{b R} h (r)^2 \dx{\tilde{P}_r} (r) \leq (n - 1) \int_0^R (n - 1) \frac{(R - r)^2}{r^2} h (r)^2 \dx{\tilde{P}_r} (r)
.\]

Combining \eqref{bound-lambda0.0}, \eqref{bound-lambda0.1} and \eqref{bound-lambda0.2}, we finally get the wanted bound \eqref{eqn:Poincare.ineq0}.
\end{proof}

\begin{acknowledgement}
The technical tools for the above proof were provided to me by Franck \textsc{Barthe}, which I warmly thank for his much precious help.
\end{acknowledgement}

\bibliographystyle{plainnat}
\bibliography{WassersteinSobolev}

\end{document}